\documentclass[11pt]{amsart}

\usepackage{amsmath}
\usepackage{mathrsfs}
\usepackage[all]{xy}
\usepackage{amssymb}
\usepackage{comment}
\usepackage{color}
\usepackage[colorlinks,citecolor=blue,urlcolor=black,linkcolor=black]{hyperref}

\newtheorem{theorem}{Theorem}[section]
\newtheorem{claim}[theorem]{Claim}

\newtheorem{lemma}[theorem]{Lemma}

\newtheorem{obs}[theorem]{Observation}

\theoremstyle{definition}
\newtheorem{definition}[theorem]{Definition}

\newtheorem{question}[theorem]{Question}

\theoremstyle{remark}

\newcount\skewfactor
\def\mathunderaccent#1#2 {\let\theaccent#1\skewfactor#2
\mathpalette\putaccentunder}
\def\putaccentunder#1#2{\oalign{$#1#2$\crcr\hidewidth
\vbox to.2ex{\hbox{$#1\skew\skewfactor\theaccent{}$}\vss}\hidewidth}}

\def\smallbox#1{\leavevmode\thinspace\hbox{\vrule\vtop{\vbox
   {\hrule\kern1pt\hbox{\vphantom{\tt/}\thinspace{\tt#1}\thinspace}}
   \kern1pt\hrule}\vrule}\thinspace}

\newcommand{\bool}{{\bf B}}

\newcommand{\cf}{{\rm cf}}

\newcommand{\Length}{{\rm Length}}

\def\qedref#1{$\qed_{\reforiginal{#1}}$}

\title{Length and ultraproducts}
\author{Shimon Garti}
\address{Institute of Mathematics
 The Hebrew University of Jerusalem,
 Jerusalem 9190401, Israel}
\email{shimon.garty@mail.huji.ac.il}

\author{Saharon Shelah}
\address{Institute of Mathematics
 The Hebrew University of Jerusalem,
 Jerusalem 9190401, Israel
 and  Department of Mathematics
 Rutgers University
 New Brunswick, NJ 08854, USA}
\email{shelah@math.huji.ac.il}
\urladdr{http://www.math.rutgers.edu/\char`\~shelah}
\thanks{Research supported by European Research Council, grant 338821. This is publication 1128 of the second author}

\subjclass[2010]{06E05, 03G05, 03E04}
\keywords{Boolean algebras, Length, pcf theory}

\begin{document}
\let\labeloriginal\label
\let\reforiginal\ref
\def\ref#1{\reforiginal{#1}}
\def\label#1{\labeloriginal{#1}}

\begin{abstract}
We shall construct, in \textsf{ZFC}, a sequence of boolean algebras $\langle\bool_i:i<\kappa\rangle$ such that $|\prod_{i<\kappa}\Length(\bool_i)/D| < \Length(\prod_{i<\kappa}\bool_i/D)$, where $D$ is a uniform ultrafilter over $\kappa$.
The method gives many instances of this inequality upon increasing the power set of a singular cardinal.
\end{abstract}

\maketitle

\newpage

\section{Introduction}

Let $\bool$ be a boolean algebra.
The length of $\bool$ is a classical cardinal invariant defined on boolean algebras and related to linearly ordered subsets of $\bool$.
The official definition reads as follows.

\begin{definition}
  \label{deflength} Length. \newline
  Let $\bool$ be a boolean algebra. The length of $\bool$ is defined by
  $\Length(\bool)=\bigcup\{|A|:A\subseteq\bool, A\ \text{is linearly-ordered by} \leq_\bool\}$.
\end{definition}

The subject of this paper is ultraproduct constructions of boolean algebras.
Let $\kappa$ be an infinite cardinal and let $D$ be a uniform ultrafilter over $\kappa$.\footnote{Uniformity means that $|A|=\kappa$ whenever $A\in{D}$.}
Given a sequence $\langle\bool_i:i<\kappa\rangle$ of boolean algebras one can form the ultraproduct algebra $\bool=\prod_{i\in\kappa}\bool_i/D$.
A classical problem (applied here to Length) is whether $|\prod_{i\in\kappa}\Length(\bool_i)/D|$ may differ from $\Length(\bool)$.
The question is investigated in general with respect to many cardinal invariants of boolean algebras, see \cite{MR3184760}.

It is consistent that $\Length(\bool)\geq|\prod_{i\in\kappa}\Length(\bool_i)/D|$ for every $\kappa$ and every uniform ultrafilter $D$ over $\kappa$.
Indeed, this inequality holds if $D$ is a regular ultrafilter, and if $V=L$ then every uniform ultrafilter is regular as proved by Donder in \cite{MR969944}.\footnote{Actually, the result of Donder holds in the core model $K$ as well.}
Thus the interesting question arises with respect to strict inequalities.

Magidor and Shelah forced in \cite{MR1674385} $\Length(\bool)<|\prod_{i\in\kappa}\Length(\bool_i)/D|$ from large cardinals in the ground model.
By the above mentioned result for regular cardinals it is clear that an inequality in this direction cannot be established in \textsf{ZFC}.
Moreover, it requires some large cardinals in the ground model.
The opposite direction, however, can be proved in \textsf{ZFC}.
A \textsf{ZFC} construction in which $\Length(\bool)>|\prod_{i\in\kappa}\Length(\bool_i)/D|$ appeared in \cite{MR1728851}.
In that paper, $\Length(\bool)=\mu^+$ where $\mu$ is a strong limit singular cardinal and $|\prod_{i\in\kappa}\Length(\bool_i)/D|=\mu$.

The construction in this paper serves as another \textsf{ZFC} example.
It is based on a different approach, rooted in \textsf{pcf} theory, and using a combinatorial property of scales over products of regular cardinals.
We indicate that the basic idea gives many more instances of the above inequality, since it enables one to get $\Length(\bool)=\lambda^+>\lambda=|\prod_{i\in\kappa}\Length(\bool_i)/D|$ for $\mu<\lambda<2^\mu$, where $\mu$ is a strong limit singular cardinal.
Of course, the additional cases do not follow from the axioms of \textsf{ZFC} since $2^\mu$ has to be larger than $\mu^+$.

The rest of the paper contains two additional sections.
In the first section we deal with ${\rm Pr}_2$, a combinatorial property that holds in \textsf{ZFC} at scales.
The main theorem of this section is interesting by its own, and we believe that it has mathematical consequences apart from length of boolean algebras.
In the second section we apply ${\rm Pr}_2$ to boolean algebras and obtain the examples related to length of ultraproducts.

Our notation is mostly standard.
We suggest \cite{MR991565} as a wonderful monograph about boolean algebras and \cite{MR3184760} for a more concrete treatment to cardinal invariants on boolean algebras.
We suggest \cite{MR2768693} as a useful source for background in pcf theory and we refer to \cite{MR1318912} for advanced material.

\newpage

\section{Pcf theory}

The main body of this section is dedicated to pcf theory.
We include here the basic definitions and facts to be used throughout the paper.
A set of cardinals $\mathfrak{a}$ is called \emph{progressive} if it consists of regular cardinals and $|\mathfrak{a}|<{\rm min}(\mathfrak{a})$.
Pcf theory analyzes the spectrum of cofinalities of the form ${\rm tcf}(\prod\mathfrak{a}/J)$ when $J$ is an ideal over $\mathfrak{a}$. By $J_{\mathfrak{a}}^{\rm bd}$ we denote the ideal of bounded subsets of $\mathfrak{a}$. The shorthand tcf stands for \emph{true cofinality}.
A central notion in pcf theory is ${\rm pcf}(\mathfrak{a})$, the collection of all true cofinalities for a specific $\mathfrak{a}$. Assume $\mathfrak{a}$ is a progressive set, and let ${\rm pcf}(\mathfrak{a})=\{\lambda:(\exists J)({\rm tcf} (\prod\mathfrak{a},<_J)=\lambda)\}$.
Occasionally we require that $\lim_J(\mathfrak{a})=\mu$, which means that $\mathfrak{a}\cap\partial\in J$ whenever $\partial<\mu$.

If $\mathfrak{a}$ is progressive then ${\rm pcf}(\mathfrak{a})$ has a last member, denoted by ${\rm max\ pcf}(\mathfrak{a})$. We call $\mathfrak{a}$ an interval if every regular cardinal between ${\rm min}(\mathfrak{a})$ and ${\rm sup}(\mathfrak{a})$ belongs to $\mathfrak{a}$. The following is a fundamental theorem of pcf theory:

\begin{theorem}
\label{ttnnoholes}
The no-holes theorem. \newline
If $\mathfrak{a}$ is a progressive interval then ${\rm pcf}(\mathfrak{a})$ is an interval as well.
\end{theorem}

\hfill \qedref{ttnnoholes}

Another central notion which we need is the pseudo power (abbreviated as pp). The name indicates that this notion may serve as a parallel to the classical notion of the power operation with respect to singular cardinals. If $\mathfrak{a}$ is a progressive set and $\kappa<|\mathfrak{a}|$ then ${\rm pcf}_\kappa(\mathfrak{a})=\bigcup\{{\rm pcf}(\mathfrak{b}): \mathfrak{b} \subseteq\mathfrak{a}\ {\rm and}\ \mathfrak{b}$ is unbounded in $\kappa\}$. If $\lambda>\cf(\lambda)=\kappa$ then ${\rm pp}_\kappa(\lambda)= \sup\{{\rm tcf}(\prod\mathfrak{a},J): \lambda=\sup(\mathfrak{a}), |\mathfrak{a}|\leq\kappa, J\supseteq J^{\rm bd}_{\mathfrak{a}}\}$.
By ${\rm pp}_\kappa^+(\lambda)$ we denote \emph{the first regular cardinal} which has no representation as the true cofinality of $\prod\mathfrak{a}/J$ for some ideal $J$ over $\mathfrak{a}$ when $|\mathfrak{a}|\leq\kappa$.

The celebrated paradox of Achilles and the tortoise comes from the school of Parmenides.\footnote{See \cite[Book VI, 239b15-30]{aristo}.}
One point is recurrent in modern treatment to this paradox, and it is connected to the separation between infinitude and boundedness.
The paradox describes an infinite sum of intervals in which Achilles tries to take over the tortoise and fails. The Greeks concluded that he will never take over the tortoise, hence the paradox.
However, infinitude need not imply this conclusion.
The sum of intervals in which Achilles cannot take over the tortoise is \emph{bounded}, despite its infinite description.
If one takes a broader look at the race beyond the bounded point of the sum of intervals then Achilles takes over the tortoise.

A progressive set is a \emph{tortoise set} if $\partial\in\mathfrak{a}$ implies ${\rm max\ pcf}(\mathfrak{a}\cap\partial) <\partial$.
The progress of Achilles is mirrored here in ${\rm pcf}(\mathfrak{a}\cap\partial)$ while ${\rm max\ pcf}(\mathfrak{a}\cap \partial)$ reflects the position of the tortoise.
The \emph{Achilles and the tortoise} theorem says that if $\mathfrak{a}$ is progressive then there is no tortoise set $\mathfrak{b}\subseteq {\rm pcf}(\mathfrak{a})$ of size $|\mathfrak{a}|^+$.
As in the original paradox, if one takes a broader look at ${\rm pcf}(\mathfrak{a})$ then eventually Achilles takes over.

A progressive set $\mathfrak{a}$ has the tortoise property when $\mathfrak{a}$ is a tortoise set. Though the size of a tortoise set within a given set of the form ${\rm pcf}(\mathfrak{a})$ is limited by $|\mathfrak{a}|^+$, one can create, in advance, a tortoise set of any size.
We shall use a theorem which ensures the tortoise property (under some assumptions on the pseudo power).

\begin{claim}
\label{ccpptortoise} pp and the tortoise property. \newline
Assume that:
\begin{enumerate}
\item [$(a)$] $\kappa<\theta=\cf(\mu)<\mu<\lambda=\cf(\lambda)$.
\item [$(b)$] ${\rm pp}_\theta^+(\mu)>\lambda$.
\item [$(c)$] $\alpha<\mu \Rightarrow \alpha^\theta<\mu$.
\item [$(d)$] $(\mu_i:i\in\theta)$ is increasing and continuous with limit $\mu$.
\end{enumerate}
Then there is a set $\mathfrak{a}\subseteq{\rm Reg}\cap\mu$ with $\min(\mathfrak{a})>{\rm otp}(\mathfrak{a})=\theta, \mu=\bigcup\mathfrak{a}$ and $\lambda={\rm tcf}(\prod\mathfrak{a},J^{\rm bd}_\theta)$ such that:
\begin{enumerate}
  \item [$(\aleph)$] If $\sigma\in\mathfrak{a}$ then $\prod(\sigma\cap\mathfrak{a})<\sigma$.
  \item [$(\beth)$] There is a scale $\bar{f}=(f_\alpha:\alpha\in\lambda)$ in the product such that $\mathcal{T}_\partial=\{f_\alpha\upharpoonright(\partial\cap\mathfrak{a}):\alpha\in\lambda\}$ has less than $\partial$-many members for every $\partial\in\mathfrak{a}$.
  \item [$(\gimel)$] There is a partition $(\mathfrak{a}_\varepsilon:\varepsilon\in\kappa)$ of $\mathfrak{a}$ into $J^{\rm bd}_\theta$-positive sets.
\end{enumerate}
\end{claim}

\par\noindent\emph{Proof}. \newline
By \cite[VII, Theorem 1.1(2)]{MR1318912} there is a set $\mathfrak{a}\subseteq{\rm Reg}\cap\mu$ as required.
We may assume that $(\aleph)$ holds upon replacing $\mathfrak{a}$ by a sub-sequence if needed.
Let $\bar{f}=(f_\alpha:\alpha\in\lambda)$ be a scale witnessing $\lambda={\rm tcf}(\prod\mathfrak{a},J^{\rm bd}_\theta)$.
The fact that $|\mathcal{T}_\partial|<\partial$ for every $\partial\in\mathfrak{a}$ follows from $(\aleph)$, thus $(\beth)$ holds true.
Part $(\gimel)$ is clear.

\hfill \qedref{ccpptortoise}

The last concept that we need is a combinatorial property called ${\rm Pr}_2$.
The following version comes from \cite{MR1318912}.
Assume that $\kappa,\theta\leq\mu\leq\lambda$.
A coloring $c:[\lambda]^2\rightarrow\kappa$ satisfies ${\rm Pr}_2(\lambda,\mu,\kappa,\theta)$ if the following holds.
Given $\xi\in\theta$ and a color $\gamma\in\kappa$, for every collection $(\alpha_{i\zeta}:\zeta\in\xi)$ of strictly increasing sequences of ordinals of $\lambda$ for each $i\in\mu$, one can find $i<j<\mu$ such that for every $\zeta_1,\zeta_2\in\xi$:
  \begin{enumerate}
    \item [$(a)$] If $\zeta_1=\zeta_2$ then $c(\alpha_{i\zeta_1},\alpha_{j\zeta_2})=\gamma$.
    \item [$(b)$] If $\zeta_1\neq\zeta_2$ then $c(\alpha_{i\zeta_1},\alpha_{j\zeta_2})=c(\alpha_{i\zeta_1},\alpha_{i\zeta_2})$.
  \end{enumerate}
Let us define a slightly different form of this concept.

\begin{definition}
  \label{defpr2} Assume that $\kappa,\sigma\leq\lambda$.
  The property ${\rm Pr}_2(\lambda,\kappa,\sigma)$ says that for some coloring $c:[\lambda]^2\rightarrow\kappa$, if $(A)$ then $(B)$, where:
  \begin{enumerate}
    \item [$(A)$] $\xi<1+\sigma$ and $\gamma_{\alpha\varepsilon}=\gamma(\alpha,\varepsilon)\in\lambda$ for every $\alpha\in\lambda,\varepsilon\in\xi$ are pairwise distinct.
    \item [$(B)$] For every $\iota\in\kappa$ there are $\alpha<\beta<\lambda$ such that:
    \begin{enumerate}
      \item [$(a)$] If $\varepsilon,\zeta<\xi$ and $\varepsilon\neq\zeta$ then $c(\gamma_{\alpha\varepsilon},\gamma_{\beta\zeta})=c(\gamma_{\alpha\zeta},\gamma_{\beta\zeta})$.
      \item [$(b)$] If $\varepsilon<\xi$ then $c(\gamma_{\alpha\varepsilon},\gamma_{\beta\varepsilon})=\iota$.
    \end{enumerate}
  \end{enumerate}
\end{definition}

Observe that ${\rm Pr}_2(\lambda,\kappa,\sigma)$ here is identical with ${\rm Pr}_2(\lambda,\lambda,\kappa,\sigma)$ in \cite{MR1318912}.
Let us add another trait to this property.

\begin{definition}
  \label{defpr2plus} Assume that $\kappa,\sigma\leq\lambda$.
  We shall say that ${\rm Pr}_2^+(\lambda,\kappa,\sigma)$ holds if ${\rm Pr}_2(\lambda,\kappa,\sigma)$ holds and, in addition, there is $S\in[\lambda]^\lambda$ such that for every distinct elements $\varepsilon,\zeta\in\xi$ the sequence $(c(\gamma_{\alpha\varepsilon},\gamma_{\beta\zeta}):\alpha,\beta\in{S},\alpha\neq\beta)$ is constant.
\end{definition}

Here is a useful fact.

\begin{obs}
  \label{obstrivial} Let $\lambda$ be a regular cardinal.
  If $\theta^{<\sigma}<\lambda$ for every $\theta\in\lambda$ then in part $(A)$ of Definition \ref{defpr2} we may assume that:
  \begin{enumerate}
    \item [$(\alpha)$] For every $\alpha<\beta$ and $\varepsilon<\zeta$ one has $\gamma_{\alpha\varepsilon}<\gamma_{\beta\zeta}$.
    \item [$(\beta)$] If $\varepsilon<\zeta$ then $\gamma_{\alpha\varepsilon}<\gamma_{\alpha\zeta}$.
  \end{enumerate}
\end{obs}

Thus we may assume that the sequences $(\gamma_{\alpha\varepsilon}:\varepsilon<\xi)$ are strictly increasing, and the elements of $(\gamma_{\beta\varepsilon}:\varepsilon<\xi)$ begin above all the elements of $(\gamma_{\alpha\varepsilon}:\varepsilon<\xi)$ whenever $\alpha<\beta$.
The latter is obtained by thinning-out the system (using the regularity of $\lambda$) and the former by renaming.

Let $(\ast)_{\lambda,\mu,\sigma,\kappa}$ be the conjunction of the following:
\begin{enumerate}
  \item [$(a)$] $\theta=\cf(\mu)<\mu<\lambda=\cf(\lambda)$.
  \item [$(b)$] $\alpha\in\lambda\Rightarrow\alpha^{<\sigma}<\lambda$, and $\sigma\leq\cf(\partial)\leq\partial$.
  \item [$(c)$] $\bar{\lambda}=(\lambda_i:i\in\partial)\subseteq{\rm Reg}\cap(\partial,\mu)$.
  \item [$(d)$] $J\supseteq J^{\rm bd}_\partial$ is a $\sigma$-complete ideal over $\partial$.
  \item [$(e)$] ${\rm tcf}(\prod_{i\in\partial}\lambda_i,J)=\lambda$, and $\kappa<\cf(\partial)$.
  \item [$(f)$] $(B_\iota:\iota\in\kappa)$ is pairwise distinct, $B_\iota\in{J^+}$ for every $\iota\in\kappa$.
  \item [$(g)$] If $j\in\partial$ then $\cf(\prod_{i\in{j}}\lambda_i)<\lambda_j$.
\end{enumerate}

\begin{claim}
\label{ccpcf} Assume $(\ast)_{\lambda,\mu,\sigma,\kappa}$.
Then there is a scale $\bar{f}=(f_\alpha:\alpha\in\lambda)$ in $(\prod_{i\in\partial}\lambda_i,J)$ such that:
\begin{enumerate}
  \item [$(\aleph)$] If $j\in\partial$ then $|\mathcal{F}_j|<\lambda_j$ where $\mathcal{F}_j=\{f_\alpha\upharpoonright{j}:\alpha\in\lambda\}$.
  \item [$(\beth)$] The function $d:[\lambda]^2\rightarrow\partial$ is defined by $d(\alpha,\beta)=\min\{j\in\partial:f_\alpha(j)\neq f_\beta(j)\}$ and $c:[\lambda]^2\rightarrow\kappa$ is defined by $c(\alpha,\beta)=\iota$ iff $d(\alpha,\beta)\in B_\iota$.
  \item [$(\gimel)$] ${\rm Pr}_2^+(\lambda,\kappa,\sigma)$ is witnessed by $c$.
\end{enumerate}
\end{claim}

\par\noindent\emph{Proof}. \newline
The existence of $\bar{f}$ satisfying $(\aleph)$ is well-known, and summarized in Claim \ref{ccpptortoise}.
Part $(\beth)$ gives the definition of $c$, so let us prove $(\gimel)$.
We are given $\xi<1+\sigma$ and $\gamma_{\alpha\varepsilon}$ is an ordinal of $\lambda$ for every $\alpha\in\lambda,\varepsilon\in\xi$.
These ordinals are pairwise distinct and by Observation \ref{obstrivial} they form a system of increasing sequences separated from each other.

For every $\alpha\in\lambda$ let $j_\alpha\in\partial$ be the first ordinal for which the sequence $(f_{\gamma(\alpha,\varepsilon)}\upharpoonright{j_\alpha}:\varepsilon<\xi)$ consists of pairwise distinct functions.
Such an ordinal exists since $\sigma\leq\cf(\partial)$.
Recall that $\sigma<\lambda=\cf(\lambda)$ and choose $j_*\in\partial, (g_\varepsilon:\varepsilon<\xi)$ and $S\in[\lambda]^\lambda$ such that $j_*=j_\alpha$ and $f_{\gamma(\alpha,\varepsilon)}\upharpoonright{j_*}=g_\varepsilon$ for every $\alpha\in{S}$.
The additional part of ${\rm Pr}_2^+(\lambda,\kappa,\sigma)$ is satisfied, as follows from the choice of $S$ and the definition of $c$.
We are left with showing that ${\rm Pr}_2(\lambda,\kappa,\sigma)$ holds.

For simplicity, assume that $S=\lambda$.
Fix $\iota\in\kappa$.
For each $j\in\partial$ let $\Lambda_j$ be the set of all sequences of the form $\bar{\beta}=(\beta_\varepsilon:\varepsilon<\xi)\in{}^\xi(\lambda_j)$ for which the set $S_{j,\bar{\beta}}=\{\alpha\in\lambda:\varepsilon<\xi\Rightarrow f_{\gamma(\alpha,\varepsilon)}(j)=\beta_\varepsilon\}$ is unbounded in $\lambda$.
Let $A$ be the set of $j\in\partial$ satisfying the following proviso:
$(\forall\beta\in\lambda_j)(\exists\bar{\beta}\in\Lambda_j)(\forall\varepsilon<\xi,\beta\leq\beta_\varepsilon)$.
Thus $A\subseteq\partial$ and we claim that $A=\partial\ \text{mod}\ J$.

To see this, let $B=\partial-A$ and assume towards contradiction that $B\in{J^+}$.
For every $j\in{B}$ choose $\beta_j\in\lambda_j$ and $\alpha_j\in\lambda$ so that if $\alpha\in[\alpha_j,\lambda)$ then $f_{\gamma(\alpha,\varepsilon}(j)<\beta_j$ for some $\varepsilon<\xi$.
Let $\gamma=\bigcup\{\alpha_j:j\in{B}\}$ and notice that $\gamma\in\lambda$ as $\partial<\lambda=\cf(\lambda)$.
Define $g\in\prod_{i\in\partial}\lambda_i$ as follows.
If $j\in{B}$ then $g(j)=\beta_j+1$ and if $j\notin{B}$ let $g(j)=0$.

Fix $\delta\in\lambda$ such that both $\gamma<\delta$ and $g<_J f_\delta$.
Notice that $\gamma_{\delta\varepsilon}\geq\delta$ for every $\varepsilon<\xi$.
For every $\varepsilon<\xi$ let $C_\varepsilon=\{i\in\partial:g(i)\geq f_{\gamma(\delta,\varepsilon)}(i)\}$, so $C_\varepsilon\in{J}$ for each $\varepsilon<\xi$ and hence $C=\bigcup\{C_\varepsilon:\varepsilon<\xi\}\in{J}$ as $J$ is $\sigma$-complete.
Fix $i\in B-C$, recalling that $B\in{J^+}$.
On the one hand, $g(i)<f_{\gamma(\delta,\varepsilon)}(i)$ for every $\varepsilon<\xi$ since $i\notin{C}$.
On the other hand, $f_{\gamma(\delta,\varepsilon)}(i)<g(i)$ for every $\varepsilon<\xi$ since $i\in{B}$, a contradiction.
We conclude, therefore, that $A=\partial\ \text{mod}\ J$.

Fix $j\in A\cap{B_\iota}$ (recall that $B_\iota\in J^+$).
Choose $\bar{\beta}_1=(\beta^1_\varepsilon:\varepsilon<\xi)\in\Lambda_j$ and then choose $\bar{\beta}_2=(\beta^2_\varepsilon:\varepsilon<\xi)\in\Lambda_j$ for which $\beta^1_\varepsilon<\beta^2_\zeta$ whenever $\varepsilon,\zeta\in\xi$.
The choice is possible since $j\in{A}$.
Finally, pick $\alpha\in S_{j,\bar{\beta}_1}$ and $\beta>\alpha$ so that $\beta\in S_{j,\bar{\beta}_2}$, this is possible since $S_{j,\bar{\beta}_2}$ is unbounded in $\lambda$.
By the above definitions, $d(\gamma_{\alpha\varepsilon},\gamma_{\beta\varepsilon})=j$ for every $\varepsilon<\xi$ and hence $c(\gamma_{\alpha\varepsilon},\gamma_{\beta\varepsilon})=\iota$ for every $\varepsilon<\xi$.
Likewise, $c(\gamma_{\alpha\varepsilon},\gamma_{\beta\zeta}) = c(\gamma_{\alpha\zeta},\gamma_{\beta\zeta})$ for every $\varepsilon,\zeta<\xi$ by the definition of $S_{j,\bar{\beta}_1}$ and $S_{j,\bar{\beta}_2}$, so we are done.

\hfill \qedref{ccpcf}

We conclude this section with a theorem of Sikorski about extending homomorphisms in boolean algebras.\footnote{A proof can be found in \cite{MR177920} and in \cite{MR991565}.}
This theorem will be used several times in our constructions.

\begin{theorem}
\label{thmsikor} If $\bool_0$ is a boolean algebra freely generated from $\{x_\gamma:\gamma\in\lambda\}$ except the inequalities mentioned in $\Gamma\subseteq\{(x_\alpha\leq x_\beta):\alpha,\beta\in\lambda\}$, $f:\{x_\gamma:\gamma\in\lambda\}\rightarrow\bool_1$ is homomorphic and $(x_\alpha\leq x_\beta)\in\Gamma\Rightarrow f(x_\alpha)\leq_{\bool_1}f(x_\beta)$ then $f$ has a homomorphic extension $\hat{f}:\bool_0\rightarrow\bool_1$.
\end{theorem}

\hfill \qedref{thmsikor}

\newpage

\section{Length}

In this section we focus on linearly ordered subsets of a boolean algebra, and we construct a sequence of boolean algebras $\langle\bool_i:i<\kappa\rangle$ such that $|\prod_{i<\kappa}\Length(\bool_i)/D| < \Length(\prod_{i<\kappa}\bool_i/D)$, where $D$ is a uniform ultrafilter over $\kappa$.
Ahead of the formal arguments, let us explain the philosophy of the proof.

There is a connection between the existence of large linearly ordered subsets of a boolean algebra and the structure of automorphisms of this algebra.
If $A\subseteq\bool$ is linearly ordered and $x,y\in A$ then there is no automorphism $f:\bool\rightarrow\bool$ such that $f(x)=y\wedge f(y)=x$.
Hence large linearly ordered subsets imply many restrictions on automorphisms and homomorphisms defined over a given boolean algebra.
Similarly, if we have a variety of automorphisms then the size of linearly ordered subsets decreases.

If one wishes to have many homomorphisms then a free boolean algebra (or any other algebraic structure) is the best way to ensure it.
If $\bool$ is a free boolean algebra then $\Length(\bool)\leq\aleph_0$ no matter how large the size of $\bool$ is.
In this case there are many homomorphisms and hence only countable linearly ordered sets.
One can control the number of homomorphisms and the size of linearly ordered sets by imposing some constraints on the freeness of the algebra.
Our algebras will be free, but we also wish to limit the possible homomorphisms, so we define a set of desired inequalities and we take algebras which are freely generated except some prescribed list of inequalities.

The purpose of these inequalities is to make sure that there are no large linearly ordered sets in each algebra.
On the other hand, we need a large linearly ordered set in the product algebra.
We define partial orders which determine the order at each algebra, and set our inequalities.
Along the definition, many of these partial orders are satisfied at every index (this will make sure that a $D$-positive set of indices will always obey these partial orders, when $D$ is a uniform ultrafilter), and hence in the product algebra there will be a large linearly ordered set.
On the other hand, at each specific index we will keep a constraint on the pertinent partial order. This will make sure that there is no very large linearly ordered set at each algebra.

More results based on the same idea can be found in \cite{MR4102862}, in which the constructions apply to both Length and Depth of boolean algebras.
Ahead of the main theorem, we define a subtler version of $\Length$.

\begin{definition}
  \label{deflengthplus} Let $\bool$ be a boolean algebra.
  The invariant $\Length^+(\bool)$ is the union of the cardinals $\theta^+$, for which there exists a linearly ordered subset $A$ of $\bool$ of size $\theta$.
\end{definition}

The advantage of $\Length^+$ is manifested in boolean algebra whose Length is a limit cardinal $\lambda$.
If $\Length(\bool)=\lambda$ and $\lambda$ is a limit cardinal (either regular or singular) then there are two different scenarios.
Either there is a linearly ordered set of size $\lambda$, in which case the Length is attained, or there is a linearly ordered set of size $\theta$ for every $\theta<\lambda$ but there is no such a set of size $\lambda$.
In both cases $\Length(\bool)=\lambda$, so the Length cannot distinguish these two situations.
However, if the Length is attained then $\Length^+(\bool)=\lambda^+$, and if the Length is not attained then $\Length^+(\bool)=\lambda$.
Thus $\Length^+$ is more accurate from this point of view.

\begin{theorem}
\label{thmmt} Assume that:
\begin{enumerate}
\item [$(\aleph)$] $\kappa<\theta=\cf(\mu)<\mu$.
\item [$(\beth)$] $\alpha<\mu\Rightarrow\alpha^\theta<\mu$.
\item [$(\gimel)$] ${\rm pp}^+_\theta(\mu)>\lambda=\cf(\lambda)$.
\end{enumerate}
Then there is a sequence $\langle\bool_i:i<\kappa\rangle$ of boolean algebras, such that:
\begin{enumerate}
\item [$(a)$] For every $i<\kappa, \Length^+(\bool_i)\leq\lambda$.
\item [$(b)$] $\Length^+(\bool)>\lambda$ where $\bool = \prod_{i<\kappa}\bool_i/D$ is the product algebra, $D$ being any uniform ultrafilter over $\kappa$.
\end{enumerate}
\end{theorem}

\par\noindent\emph{Proof}. \newline
Since $\lambda<{\rm pp}^+_\theta(\mu)$ one can find a sequence $\bar{\lambda} = \langle\lambda_i:i<\theta\rangle$ and an ideal $I$ over $\theta$ such that:
\begin{enumerate}
\item [$(a)$] $\lambda_i=\cf(\lambda_i)<\mu$ for every $i<\theta$.
\item [$(b)$] $\lim_I(\bar{\lambda})=\mu$.
\item [$(c)$] ${\rm tcf}(\prod_{i<\theta}\lambda_i,I)\geq\lambda$.
\end{enumerate}
We may assume that ${\rm tcf}(\prod_{i<\theta}\lambda_i,I)=\lambda$ by virtue of Theorem \ref{ttnnoholes}.
Without loss of generality, $\bar{\lambda}$ is increasing, $I=J^{\rm bd}_\theta$ and ${\rm max\ pcf}\{\lambda_j:j<i\}<\lambda_i$ for every $i<\theta$.

Let $\mathfrak{a} = \{\lambda_i:i<\theta\}$, and decompose $\mathfrak{a}$ into $(\mathfrak{a}_i:i<\kappa)$ such that each $\mathfrak{a}_i$ is $I$-positive and $i<\kappa\Rightarrow\{j\in\theta:\lambda_j\in\mathfrak{a}_i\}$ is a stationary subset of $\theta$.
Define $J = \{\mathfrak{b}\subseteq\mathfrak{a}: \{j<\theta:\lambda_j\in \mathfrak{b}\}\in I\}$.
So $J$ is an ideal over $\mathfrak{a}$ and $J\supseteq J^{\rm bd}_{\mathfrak{a}}$.
Further, ${\rm tcf}(\prod_{j<\theta}\lambda_j,I) = {\rm tcf}(\prod \mathfrak{a}, J) = \lambda$.
Using Claim \ref{ccpptortoise} we choose a sequence $\bar{f} = (f_\alpha: \alpha\in\lambda)\subseteq\prod\mathfrak{a}$, increasing and cofinal in $(\prod\mathfrak{a},J)$, such that $|\mathcal{T}_\partial|<\partial$ for every $\partial\in\mathfrak{a}$, where $\mathcal{T}_\partial = \{f_\alpha\upharpoonright(\mathfrak{a}\cap\partial): \alpha\in\lambda\}$.

Fix any uniform ultrafilter $D$ over $\kappa$.
For $\alpha<\beta<\lambda$ let $\partial_{\alpha\beta} = \min\{\partial\in\mathfrak{a}: f_\alpha(\partial)\neq f_\beta(\partial)\}$.
We define a binary relation $<_*$ over $\lambda$ by $\alpha<_*\beta$ iff $f_\alpha(\partial_{\alpha\beta})<f_\beta(\partial_{\alpha\beta})$.
Based on $<_*$ we define a binary relation $<_i$ over $\lambda+1\cup\{-1\}$, for each $i<\kappa$, as follows:
\begin{gather*}
\alpha<_i\beta \Leftrightarrow \\
(\alpha<_*\beta \wedge \partial_{\alpha\beta}\in \bigcup \{\mathfrak{a}_j:j<i\}) \vee \\
(\alpha=-1 \wedge \beta\in\lambda+1) \vee \\
(\alpha\in\lambda\cup\{-1\} \wedge \beta=\lambda)
\end{gather*}
It is routine to verify that $<_*$ is a linear order and each $<_i$ is a partial order. Likewise, if $\alpha<_*\beta$ then $\{i\in\kappa:\alpha<_i\beta\}\in D$, since it is an end-segment of $\kappa$ and $D$ is uniform.

Let $\bool_i$ be the boolean algebra generated freely from $\{x_\alpha:\alpha\in\lambda+1\cup\{-1\}\}$ except the inequalities mentioned in $\Gamma_i$, where $\Gamma_i = \{(x_\alpha\leq x_\beta):\alpha\leq_i\beta\} \cup \{x_{-1}=0,x_\lambda=1\}$.
Put another way, $\bool_i\models x_\alpha<x_\beta \Leftrightarrow \alpha<_i\beta$.
Let $\bool$ be the product algebra $\prod_{i<\kappa}\bool_i/D$.

It is easy to see that $\Length^+(\bool)>\lambda$.
For this, we point to a linearly ordered subset of $\bool$ of size $\lambda$.
For every $\alpha\in\lambda$ let $c_\alpha\in\prod_{i<\kappa}\bool_i$ be the constant function with value $x_\alpha$, and let $c_\alpha/D\in\bool$ be its $D$-equivalence class.
Clearly $\alpha\neq\beta\Rightarrow c_\alpha/D\neq c_\beta/D$, hence the set $A = \{c_\alpha/D:\alpha\in\lambda\}\subseteq\bool$ is of size $\lambda$.
Now if $\alpha<_*\beta$ then $\{i<\kappa:\bool_i\models x_\alpha<x_\beta\} = \{i<\kappa: \alpha<_i\beta\}\in D$, being an end-segment of $\kappa$.
Consequently, $\alpha<_*\beta\Rightarrow c_\alpha/D<_{\bool}c_\beta/D$.
Since $<_*$ is a linear order we see that $A$ is a linearly ordered subset of $\bool$ and hence $\Length(\bool)\geq\lambda$, so $\Length^+(\bool)>\lambda$.

The other side of the coin is that $\Length^+(\bool_i)\leq\lambda$ for every $i<\kappa$.
For proving this, fix an ordinal $i\in\kappa$ and assume toward a contradiction that $\lambda<\Length^+(\bool_i)$ for some $i<\kappa$.
Fix a linearly ordered set $\{a_\alpha:\alpha<\lambda\}\subseteq\bool_i$ which exemplifies this statement.

Each $a_\alpha$ is an element of the boolean algebra $\bool_i$ and hence expressible as $\sigma_\alpha(x_{\gamma(\alpha,0)},\ldots,x_{\gamma(\alpha,n(\alpha)-1)})$, where $\sigma_\alpha$ is a boolean term and $\gamma(\alpha,\ell)\in\lambda+1\cup\{-1\}$.
But $\lambda=\cf(\lambda)>\aleph_0$, and the number of boolean terms is $\aleph_0$.
So without loss of generality $\sigma_\alpha=\sigma$ and $n(\alpha)=n$ for every $\alpha\in\lambda$, where $\sigma$ is a fixed boolean term and $n$ is a fixed natural number.
We can write $a_\alpha = \sigma(x_{\gamma(\alpha,0)},\ldots,x_{\gamma(\alpha,n-1)})$ for every $\alpha\in\lambda$, and we may assume that $\langle \gamma(\alpha,\ell):\ell<n\rangle$ is always an increasing sequence of ordinals.

For every $\alpha\in\lambda$ there is a set $p_\alpha\subseteq n\times n$ such that $(k,\ell)\in p_\alpha$ iff $\gamma(\alpha,k)<_i\gamma(\alpha,\ell)$.
Being a subset of $n\times n$, it is a finite set.
Hence for $\lambda$ many ordinals we have $p_\alpha = p$ for some fixed $p\subseteq n\times n$.
By concentrating on these ordinals we may assume that $\alpha<\beta<\lambda \Rightarrow p_\alpha=p_\beta$.
This means that $\gamma(\alpha,k)<_i\gamma(\alpha,\ell)$ iff $\gamma(\beta,k)<_i\gamma(\beta,\ell)$ for every $\alpha,\beta\in\lambda$.

Since $\bar{f}$ is an increasing sequence in $(\prod\mathfrak{a},J)$, for every $\alpha\in\lambda$ there exists $\partial_\alpha\in\mathfrak{a}$ such that the elements of $\langle f_{\gamma(\alpha,\ell)}\upharpoonright (\mathfrak{a}\cap\partial_\alpha):\ell<n\rangle$ are pairwise distinct.
Since $|\mathfrak{a}|=\theta<\cf(\lambda)=\lambda$ we may assume without loss of generality that $\partial_\alpha=\partial$ for some fixed $\partial\in\mathfrak{a}$ and each $\alpha\in\lambda$.

We apply the Delta-system lemma to the collection of $\lambda$ finite sequences $\langle\gamma(\alpha,\ell):\ell<n\rangle$. We can assume now that there exists $m<n$ such that $\ell<m\Rightarrow \gamma(\alpha,\ell)=\gamma(\ell)$ and if $\alpha<\beta$ then $\gamma(\alpha,n-1)<\gamma(\beta,m)$.
Recall that $i\in\kappa$ is fixed, and choose $j\in(i,\kappa)$ bearing in mind that $\mathfrak{a}_j\in J^+$.
From Theorem \ref{ccpcf}, one can find $\alpha<\beta<\lambda$ such that for every $k,\ell\in[m,n)$ it is true that $\partial_{\gamma(\alpha,k)\gamma(\beta,\ell)}\in\mathfrak{a}_j$.
We include the root $\langle\gamma(\alpha,\ell):\ell<m\rangle$ of the Delta-system in the first sequence, so the property ${\rm Pr}_2$ holds with respect to these elements (and all the other disjoint sequences) as well.

Define $S=\{\gamma(\alpha,\ell),\gamma(\beta,\ell):\ell<n\} \cup \{-1,\lambda\}$, so $S$ is a finite subset of $\lambda+1\cup\{-1\}$.
Observe that $(S,<_i)$ is a finite partial order with $-1,\lambda\in S$ and $s\in S\Rightarrow -1\leq_i s\leq_i \lambda$.
Hence any partial ordering $P\supseteq(S,<_i)$ projects onto $(S,<_i)$.
In particular, there exists a function $h:\lambda+1\cup\{-1\}\rightarrow S$ such that $h\upharpoonright S={\rm id}_S$ and $\gamma_0\leq_i\gamma_1 \Rightarrow h(\gamma_0)\leq_i h(\gamma_1)$.

Let $\bool_S$ be the (finite) boolean algebra generated freely from $\{x_\gamma:\gamma\in S\}$ except the relations in $\Gamma_S = \{(x_\gamma\leq x_\delta):\gamma\leq_i\delta, \gamma,\delta\in S\}$.
Define a mapping $f:\{x_\gamma:\gamma\in\lambda\}\rightarrow\bool_S$ by $f(x_\gamma)=x_{h(\gamma)}$.
By Theorem \ref{thmsikor} there is a homomorphism $\hat{f}:\bool_i\rightarrow\bool_S$ such that $f\subseteq\hat{f}$ and hence $\hat{f}(x_\gamma)=x_\gamma$ whenever $\gamma\in S$.

Define a function $g:S\rightarrow S$ by letting, for every $\ell<n$, $g(\gamma(\alpha,\ell))=\gamma(\beta,\ell)$ and $g(\gamma(\beta,\ell))=\gamma(\alpha,\ell)$.
Notice that $g$ is a permutation of $S, g\circ g={\rm id}_S$ and $g$ maps $\Gamma_S$ onto itself, namely $(x_\gamma\leq x_\delta)\in\Gamma_S\Rightarrow (x_{g(\gamma)}\leq x_{g(\delta)})\in\Gamma_S$.
Indeed, if $k,\ell\in[m,n)$ then $\gamma(\alpha,k)$ and $\gamma(\beta,\ell)$ are $\leq_i$-incomparable, so if $\sigma=\gamma(\alpha,k)$ and $\tau=\gamma(\beta,\ell)$ then $x_\sigma\leq x_\tau\notin\Gamma_S$.
If $k<m$ or $\ell<m$ then $\partial_{\gamma(\alpha,k)\gamma(\alpha,\ell)} = \partial_{\gamma(\alpha,k)\gamma(\beta,\ell)}$ by ${\rm Pr}_2(\lambda,\aleph_0,\theta)$.
Hence $\gamma(\alpha,k)<_i\gamma(\alpha,\ell)$ iff $\gamma(\alpha,k)<_i\gamma(\beta,\ell)$.
Similarly, $\gamma(\beta,k)<_i\gamma(\alpha,\ell)$ iff $\gamma(\beta,k)<_i\gamma(\beta,\ell)$.
Finally, $\gamma(\alpha,k)<_i\gamma(\alpha,\ell)$ iff $\gamma(\beta,k)<_i\gamma(\beta,\ell)$ since $p_\alpha=p_\beta=p$, so $\gamma(\alpha,k)<_i\gamma(\beta,\ell)$ iff $\gamma(\beta,k)<_i\gamma(\alpha,\ell)$.
Together, $x_\sigma\leq x_\tau$ iff $x_{g(\sigma)}\leq x_{g(\tau)}$, as sought.

By another application of Theorem \ref{thmsikor} there exists an automorphism $\hat{g}:\bool_S\rightarrow\bool_S$ which satisfies $\hat{g}(x_\gamma)=x_{g(\gamma)}$.
Observe that $\hat{g}$ respects boolean terms in the following sense:
\begin{center}
$\hat{g}(\sigma(\ldots,x_{\gamma(\alpha,\ell)},\ldots)) = \sigma(\ldots,x_{\gamma(\beta,\ell)},\ldots)$. \\
$\hat{g}(\sigma(\ldots,x_{\gamma(\beta,\ell)},\ldots)) = \sigma(\ldots,x_{\gamma(\alpha,\ell)},\ldots)$.
\end{center}

Back to the chosen pair of ordinals $\alpha<\beta<\lambda$, since $\{a_\gamma:\gamma\in\lambda\}$ is $<_{\bool_i}$-linearly ordered we know that $(a_\alpha<_{\bool_i}a_\beta)\vee(a_\beta<_{\bool_i}a_\alpha)$.
Without loss of generality $a_\alpha<_{\bool_i}a_\beta$, so $\bool_i \models \sigma(\ldots,x_{\gamma(\alpha,\ell)},\ldots) < \sigma(\ldots,x_{\gamma(\beta,\ell)},\ldots)$.
Since $\hat{f}:\bool_i\rightarrow\bool_S$ is homomorphic, we see that:
$$\bool_S \models \sigma(\ldots,x_{\gamma(\alpha,\ell)},\ldots) < \sigma(\ldots,x_{\gamma(\beta,\ell)},\ldots).$$
Likewise, $\hat{g}:\bool_S\rightarrow\bool_S$ is homomorphic and hence we also have:
$$\bool_S \models \sigma(\ldots,x_{\gamma(\beta,\ell)},\ldots) < \sigma(\ldots,x_{\gamma(\alpha,\ell)},\ldots).$$
This contradiction accomplishes the proof of the theorem.

\hfill \qedref{thmmt}

The above theorem deals with $\Length^+$, but we can deduce the following result about $\Length$:

\begin{theorem}
\label{thmmt1} Assume that:
\begin{enumerate}
\item [$(\aleph)$] $\kappa<\theta=\cf(\mu)<\mu$.
\item [$(\beth)$] $\alpha<\mu\Rightarrow\alpha^\theta<\mu$.
\item [$(\gimel)$] $\lambda=\mu^+$.
\end{enumerate}
Then there exists a sequence of boolean algebras $\langle\bool_i: i<\kappa\rangle$ such that:
\begin{enumerate}
\item [$(a)$] $|\bool_i|=\lambda$ for every $i<\kappa$.
\item [$(b)$] $\Length(\bool_i)\leq\mu$ for every $i<\kappa$.
\item [$(c)$] If $D$ is a uniform ultrafilter over $\kappa$ and $\bool = \prod_{i<\kappa}\bool_i/D$ is the product algebra then $\Length(\bool)\geq\lambda$.
\item [$(d)$] $\prod_{i<\kappa}\Length(\bool_i)/D=\mu$, so $\prod_{i<\kappa}\Length(\bool_i)/D < \Length(\bool)$.
\end{enumerate}
\end{theorem}

\par\noindent\emph{Proof}. \newline
Since $\lambda=\mu^+, {\rm pp}^+_\theta(\mu)>\lambda$ and hence Theorem \ref{thmmt} applies.
Let $\langle\bool_i: i<\kappa\rangle$ be a sequence of boolean algebras as asserted in Theorem \ref{thmmt}.
For every $i<\kappa$ we have $\Length^+(\bool_i)\leq\lambda$, hence (recalling that $\lambda=\mu^+$) $\Length(\bool_i)\leq\mu$.
On the other hand, $\Length^+(\bool)>\lambda$ so $\Length(\bool)\geq\lambda$. It follows that $\prod_{i<\kappa}\Length(\bool_i)/D\leq\mu^\kappa = \mu$ (the equality $\mu^\kappa=\mu$ is due to $(\beth)$), and $\mu<\lambda=\Length(\bool)$ so the proof is accomplished.

\hfill \qedref{thmmt1}

A detailed examination of the construction yields accurate values of Length for $\bool$ and each $\bool_i$ in the above theorems.
This is expressed by the following two lemmata.

\begin{lemma}
\label{lemb} In the construction of Theorem \ref{thmmt}, if $\lambda=\cf(\lambda)<\mu^{+\omega}$ then $\Length(\bool)=\lambda$.
\end{lemma}

\par\noindent\emph{Proof}. \newline
We have seen already that the set $A = \{c_\alpha/D:\alpha\in\lambda\} \subseteq \bool$ is of size $\lambda$ and it is linearly ordered by $<_\bool$, so we have $\Length(\bool)\geq\lambda$.
We claim that if $\lambda<\mu^{+\omega}$ then $|\bool|=\lambda$ and hence $\Length(\bool)\leq\lambda$ as well.
For this notice that $\mu^\kappa=\mu$ and hence by induction on $\lambda\in(\mu,\mu^{+\omega})$ we have $\lambda^\kappa=\lambda$.
Since $|\bool_i|=\lambda$ for every $i<\kappa$ we see that $\lambda\leq|\bool|=|\prod_{i<\kappa}\bool_i/D|\leq\lambda^\kappa=\lambda$, so we are done.

\hfill \qedref{lemb}

More generally one can prove that if $\lambda\rightarrow(\lambda)^2_{2^\kappa}$ then $\Length(\bool)\leq\lambda$ even if $\lambda>\mu^{+\omega}$.

\begin{lemma}
\label{lembi} In the construction of Theorem \ref{thmmt1}, $\Length(\bool_i)=\mu$ for every $i<\kappa$.
\end{lemma}

\par\noindent\emph{Proof}. \newline
Fix an ordinal $i\in\kappa$ and define the following set:
$$
A = \{\partial\in\mathfrak{a}:\sup\{f_\alpha(\partial):\alpha\in\lambda\} = \partial\}.
$$
We may assume that $A=\mathfrak{a}$ mod $J$, as this requirement can be added to the choice of $\bar{f}$.
In particular, $\sup(A)=\mu$.

Fix any $\partial\in A$ such that $\partial\in\bigcup\{\mathfrak{a}_j:j<i\}$, and notice that these elements are unbounded in $\mu$.
Choose a sequence of ordinals $\langle\alpha_\varepsilon:\varepsilon\in\partial\rangle$ so that $\langle f_{\alpha_\varepsilon}(\partial):\varepsilon\in\partial\rangle$ is an increasing sequence of ordinals in $\partial$.
Recall that $|\{f_\alpha\upharpoonright(\mathfrak{a}\cap\partial): \alpha<\lambda\}|\leq$ max pcf$(\mathfrak{a}\cap\partial)<\partial$, and $\partial=\cf(\partial)$.
Hence without loss of generality there is some fixed $g\in\prod(\mathfrak{a}\cap\partial)$ such that $\varepsilon<\partial\Rightarrow f_{\alpha_\varepsilon}\upharpoonright(\mathfrak{a}\cap\partial)=g$.

Now if $\varepsilon<\zeta<\partial$ and $\partial\in\bigcup\{\mathfrak{a}_j:j<i\}$ then, by the choice of $g$, $\partial$ is the first point in which $f_{\alpha_\varepsilon}(\partial) \neq f_{\alpha_\zeta}(\partial)$.
Moreover, $f_{\alpha_\varepsilon}(\partial)<f_{\alpha_\zeta}(\partial)$ by the choice of $\langle f_{\alpha_\varepsilon}(\partial):\varepsilon\in\partial\rangle$ and hence $\alpha_\varepsilon<_i\alpha_\zeta$.
By the definition of the order in $\bool_i$ we see that $x_{\alpha_\varepsilon}<_{\bool_i}x_{\alpha_\zeta}$ so $\langle x_{\alpha_\varepsilon}:\varepsilon\in\partial\rangle$ implies that $\Length(\bool_i)\geq\partial$.
Applying this argument to every $\partial\in{A}$ one concludes that $\Length(\bool_i)\geq\mu$, since $A$ is unbounded in $\mu$.
On the other hand, in Theorem \ref{thmmt1} we have seen that $\Length^+(\bool_i)\leq\lambda=\mu^+$.
It follows that any linearly ordered subset of $\bool_i$ is of size at most $\mu$, so we are done.

\hfill \qedref{lembi}

The case in which $\kappa<\theta=\cf(\mu)<\mu$ and $\lambda=\mu^+$ yields examples of $|\prod_{i<\kappa}\Length(\bool_i)/D| < \Length(\prod_{i<\kappa}\bool_i/D)$ in \textsf{ZFC}.
Indeed, if $\mu$ is a strong limit cardinal then $\mu^\kappa=\mu$ (actually, $\alpha^\kappa<\mu$ whenever $\alpha\in\mu$ is sufficient) and scales of length $\mu^+$ are available in \textsf{ZFC}.
We indicate, however, that if $\lambda>\mu^+$ and not so large (e.g., $\lambda<\mu^{+\kappa}$) then similar examples are at hand.
Thus if $\lambda=\chi^+$ and $\chi^\kappa=\chi$ then the main result of this section provides a sequence $(\bool_i:i\in\kappa)$ so that $\chi^\kappa=\chi=|\prod_{i<\kappa}\Length(\bool_i)/D|$ while $\Length(\prod_{i<\kappa}\bool_i/D)=\chi^+$.
Anyways, we know to build examples with a gap of one cardinality only.

\begin{question}
  \label{qlargegap} Is it consistent that there is an infinite cardinal $\kappa$ and a sequence of boolean algebras $(\bool_i:i\in\kappa)$ such that $|\prod_{i<\kappa}\Length(\bool_i)/D| < \Upsilon < \Length(\prod_{i<\kappa}\bool_i/D)$ for some infinite cardinal $\Upsilon$?
\end{question}

Note that if our coloring $c$ satisfies ${\rm Pr}_2(\lambda,\mu,\kappa,\theta)$ when $\mu>\cf(\mu)=\theta>\kappa$ is a strong limit cardinal and ${\rm pp}^+(\mu)>\lambda=\cf(\lambda)>\mu^+$ then we will get a positive answer to the above question (with $\Upsilon=\mu^+$).
So far, however, we know how to get ${\rm Pr}_2(\lambda,\lambda,\kappa,\theta)$ but not the apparently stronger statement ${\rm Pr}_2(\lambda,\mu,\kappa,\theta)$.

\newpage

\bibliographystyle{alpha}
\bibliography{arlist}

\end{document}